\documentclass[12pt]{article}
\usepackage{mathrsfs}
\usepackage{amsmath,amsfonts,amssymb,rotating}

\def\qed{\nopagebreak\hfill{\rule{4pt}{7pt}}}
\def\proof{\noindent {\it{Proof.} \hskip 2pt}}

\newtheorem{theo}{Theorem}[section]

\newtheorem{lemm}[theo]{Lemma}

\newtheorem{coro}[theo]{Corollary}

\newdimen\Squaresize \Squaresize=11pt
\newdimen\Thickness \Thickness=0.7pt
\def\Square#1{\hbox{\vrule width \Thickness
   \vbox to \Squaresize{\hrule height \Thickness\vss
    \hbox to \Squaresize{\hss#1\hss}
   \vss\hrule height\Thickness}
\unskip\vrule width \Thickness} \kern-\Thickness}

\def\Vsquare#1{\vbox{\Square{$#1$}}\kern-\Thickness}

\def\moins{\raise 1pt\hbox{{$\scriptstyle -$}}}

\begin{document}

\begin{center}
{\large \bf Recurrence Relations for

Strongly $q$-Log-Convex Polynomials}
\end{center}

\begin{center}
William Y. C. Chen$^{1}$, Larry X.W. Wang$^{2}$ and Arthur
L. B. Yang$^{3}$\\[6pt]
Center for Combinatorics, LPMC-TJKLC\\
Nankai University, Tianjin 300071, P. R. China\\
Email: $^{1}${\tt chen@nankai.edu.cn}, $^{2}${\tt
wxw@cfc.nankai.edu.cn}, $^{3}${\tt yang@nankai.edu.cn}
\end{center}

\vspace{0.3cm} \noindent{\bf Abstract.} We consider a class of
strongly $q$-log-convex polynomials based on a triangular recurrence
relation with linear coefficients, and we show that the Bell
polynomials, the Bessel polynomials, the Ramanujan polynomials and
the Dowling polynomials are strongly $q$-log-convex. We also prove
that the Bessel transformation preserves log-convexity.





\section{Introduction}

The main objective of this paper is to study the $q$-log-convexity
of a class of polynomials  whose coefficients satisfy a triangular
recurrence relation with linear coefficients. The notion of
log-convexity is closely related to log-concavity. Stanley
introduced the concept of $q$-log-concavity, which naturally leads
to the notion of $q$-log-convexity. Compared with $q$-log-concave
polynomials, $q$-log-convex polynomials did not draw much attention.
Only till recently, Liu and Wang  \cite{liuwan2006} have shown that
some classical polynomials are $q$-log-convex, such as the Bell
polynomials and the Eulerian polynomials.

In this paper, we will show that a polynomial
\begin{equation}\label{pnx}
 P_n(q)=\sum_{k=0}^n T(n,k)q^k
 \end{equation}
  is
$q$-log-convex if the coefficients $T(n,k)$ satisfy certain
recurrence relation with linear coefficients in $n$ and $k$. The
concept of strong $q$-log-concavity is due to Sagan
\cite{sagan1992t}. In this framework, we will show that the Bell
polynomials, the Bessel polynomials, the Ramanujan polynomials and
the Dowling polynomials are strongly $q$-log-convex.

Let us give a brief review of the background and terminology.
Unimodal and log-concave sequences and polynomials often arise in
combinatorics, algebra and geometry; See the surveys of Brenti
\cite{brenti1989,brenti1994} and Stanley \cite{stanley1989}. A
sequence $\{z_k\}_{k\geq 0}$ of nonnegative real numbers is said to
be unimodal if there exists an integer $r\geq 0$ such that
$$z_0\leq z_1\leq\cdots\leq z_r\geq z_{r+1}\geq
z_{r+2} \geq \cdots.$$ It is said to be log-concave if
$$z_m^2\geq z_{m+1}z_{m-1},\quad  m\geq 1,$$
and it is said to be strongly log-concave if
$$z_m z_n\geq z_{m-1}z_{n+1},\quad n\geq m\geq 1.$$ For a sequence of positive real numbers,
log-concavity is equivalent to  strong log-concavity and implies
unimodality.

Analogously, a sequence $\{z_k\}_{k\geq 0}$ of nonnegative real
numbers is said to be log-convex if
$$z_m^2\leq z_{m+1}z_{m-1}, \quad  m\geq 1,$$
and it is said to be strongly log-convex if
$$z_m z_n\leq z_{m-1}z_{n+1},\quad  n\geq m\geq 1.$$
It is also easily seen that, for a sequence of positive real
numbers, log-concavity is equivalent to the log-convexity of the
sequence of the reciprocals; See \cite{liuwan2006}.

However, the equivalence of log-concavity and strong log-concavity
does not apply to polynomial sequences. The $q$-log concavity of
polynomials has been extensively studied; See, for example,  Butler
\cite{butler1990}, Krattenthaler \cite{kratte1989}, Leroux
\cite{leroux1990}, and Sagan \cite{sagan1992t}. Adopting the
notation of Sagan \cite{sagan1992t}, we write
$$f(q)\geq_q g(q)$$ if the difference $f(q)-g(q)$ has nonnegative
coefficients as a polynomial of $q$. A sequence of polynomials
$\{f_k(q)\}_{k\geq 0}$ over the field of real numbers is called
$q$-log-concave if
$$f_{m}(q)^2\geq_q f_{m+1}(q)f_{m-1}(q),\quad  m\geq 1,$$
and it is strongly $q$-log-concave if
$$f_{m}(q)f_{n}(q)\geq_q f_{m-1}(q)f_{n+1}(q),\quad  n\geq m\geq 1.$$

Based on the $q$-log concavity, it is natural to define the
$q$-log-convexity and the strong $q$-log-convexity. We say that a
polynomial sequence $\{f_n(q)\}_{n\geq 0}$ is $q$-log-convex if
$$f_{m+1}(q)f_{m-1}(q)\geq_q f_{m}(q)^2,\quad m\geq
1,$$ and it is strongly $q$-log-convex if
$$f_{m-1}(q)f_{n+1}(q)\geq_q f_{m}(q)f_{n}(q),\quad  n\geq m\geq 1.$$
For a sequence of polynomials, $q$-log-convexity is also not
equivalent to strong $q$-log-convexity. Note that Butler and
Flanigan \cite{butfla2007} defined a different  $q$-analogue of
log-convexity.

For a  $q$-log-convex sequence of polynomials $P_n(q)$ as given in
(\ref{pnx}), we will be concerned with the linear transformation
associated with $P_n(q)$ which transforms a sequence $\{z_n\}_{n\geq
0}$ into a sequence $\{w_n\}_{n\geq 0}$ given by
$$w_n=\sum_{k=0}^nT(n,k)z_k.$$ We
say that the linear transformation preserves log-convexity, if, for
any given log-convex sequence $\{z_n\}_{n\geq 0}$ of positive real
numbers, the sequence $\{w_n\}_{n\geq 0}$ defined by the above
transformation is also log-convex. For the Bell polynomials, the
corresponding linear transformation is defined by the Stirling
numbers of the second kind. It has been proved that the Bell
transformation preserves log-convexity \cite{liuwan2006}. In this
paper, we will show that the Bessel transformation  preserves
log-convexity.

\section{The strong $q$-log-convexity}

In this section,  we consider polynomials
$$P_n(q)=\sum_{k=0}^n T(n,k)q^k,\quad n\geq 0,$$ where the
coefficients $T(n,k)$ are nonnegative real numbers and satisfy the
following recurrence relation
\begin{eqnarray}
T(n,k)&= &(a_1n+a_2k+a_3)T(n-1,k) \nonumber \\[6pt]
        & & \quad +\,(b_1n+b_2k+b_3)T(n-1,k-1),  \mbox{ for } n\geq k\geq 1, \label{rec}
\end{eqnarray}
and  the boundary conditions
$$T(n,-1)=T(n,n+1)=0, \mbox{ for } n\geq 1,$$
$$a_1\geq 0,\quad a_1+a_2\geq 0,\quad  a_1+a_2+a_3> 0,
$$
and
$$
b_1\geq 0, \quad b_1+b_2\geq 0,\quad b_1+b_2+b_3> 0.
$$
For the triangular array $\{T(n,k)\}_{n\geq k\geq 0}$,  we always
assume that $T(0,0)>0$. Thus we have $T(n,k)>0$ for $0\leq k \leq
n$.

The following lemma is a special case of  Theorem 2 of Kurtz
\cite{jcta1972}.

\begin{lemm}\label{lca}
Suppose that the positive array $\{T(n,k)\}_{n\geq k\geq 0}$
satisfies the recurrence relation \eqref{rec}. Then, for given $n$,
the sequence $\{T(n,k)\}_{0\leq k\leq n}$ is log-concave, namely,
for $0\leq k \leq n$,
\begin{equation}\label{tnk2}
T(n,k)^2\geq T(n,k-1)T(n,k+1).
\end{equation}
\end{lemm}

Using the log-concavity  (\ref{tnk2}) for the triangular array
$\{T(n,k)\}_{n\geq k\geq 0}$, Liu and Wang obtained a sufficient
condition  for the polynomial sequence $\{P_n(q)\}_{n\geq 0}$ to be
$q$-log-convex \cite[Theorem 4.1]{liuwan2006}.

\begin{theo} \label{liuwangthm} Suppose that the
array $\{T(n,k)\}_{n\geq k\geq 0}$ of positive  numbers satisfies
the recurrence relation \eqref{rec} and the additional condition
\[
(a_2b_1-a_1b_2)n+a_2b_2k+(a_2b_3-a_3b_2)\geq 0, \quad \mbox{for
$0<k\leq n$}.
\]
Then the polynomials $P_n(q)$ form a $q$-log-convex sequence.
\end{theo}

This theorem can be used to show that the Bell polynomials and the
Eulerian polynomials are $q$-log-convex. We will give alternative
conditions for the recurrence relation \eqref{rec} and will show
that our conditions are satisfied by the Bell polynomials, the
Bessel polynomials, the Ramanujan polynomials and the Dowling
polynomials. It is also hoped that further studies will be carried
out for more other types of recurrence relations with polynomial or
even rational coefficients in $n$ and $k$.

An important property of the triangular array satisfying our
conditions is described by the following lemma.

\begin{lemm}\label{change} Suppose that the array $\{T(n,k)\}_{n\geq k\geq 0}$
of positive numbers  satisfies \eqref{rec} with $a_2,b_2\geq 0$.
Then, for any $l'\geq l\geq 0$ and $m'\geq m\geq 0$, we have
\begin{equation}\label{eqch}
T(m,l)T(m',l')-T(m,l')T(m',l)\geq 0.
\end{equation}
\end{lemm}

\proof Restate (\ref{eqch}) as
\[
\frac{T(m,l')}{T(m,l)}\leq \frac{T(m',l')}{T(m',l)}.
\]
It suffices to show that for any $s\geq m$
\[
\frac{T(s,l')}{T(s,l)}\leq \frac{T(s+1,l')}{T(s+1,l)}.
\]
Let
$$f(l,l',s)=T(s,l)T(s+1,l')-T(s,l')T(s+1,l).$$
From the recurrence relation \eqref{rec}, we see that
$$
\begin{array}{rcl}
f(l,l',s)&=& \left(a_1(s+1)+a_2l'+a_3\right)T(s,l')T(s,l)\\[5pt]
          & &\quad  +\,\left(b_1(s+1)+b_2l'+b_3\right)T(s,l'-1)T(s,l)\\[8pt]
          & & \quad -\,\left(a_1(s+1)+a_2l+a_3\right)T(s,l)T(s,l')\\[8pt]
          & & \quad -\,\left(b_1(s+1)+b_2l+b_3\right)T(s,l-1)T(s,l')\\[10pt]
          &=& \left(b_1(s+1)+b_2l'+b_3\right)T(s,l'-1)T(s,l)\\[8pt]
          & & \quad -\,\left(b_1(s+1)+b_2l+b_3\right)T(s,l-1)T(s,l')\\[8pt]
          & & \quad +\,a_2(l'-l)T(s,l)T(s,l')\\[10pt]
          &\geq & \left(b_1(s+1)+b_2l'+b_3\right)T(s,l-1)T(s,l') \\[8pt]
          & & \quad -\,\left(b_1(s+1)+b_2l+b_3\right)T(s,l-1)T(s,l')\\[8pt]
          & & \quad +\,a_2(l'-l)T(s,l)T(s,l') \hfill{\rule{88pt}{0pt}(\mbox{by Lemma \ref{lca}})}\\[10pt]
          &=& b_2(l'-l)T(s,l-1)T(s,l')+a_2(l'-l)T(s,l)T(s,l'),
\end{array}
$$
Which is nonnegative in view of the condition $a_2,b_2\geq 0$. This
completes the proof. \qed

The main result of this paper is given below.

\begin{theo}\label{strong}
Suppose that the array $\{T(n,k)\}_{n\geq k\geq 0}$  of positive
numbers satisfies \eqref{rec} with  $a_2,b_2\geq 0$. Then the
polynomial sequence $\{P_n(q)\}_{n\geq 0}$ is strongly
$q$-log-convex,  namely, for any $n\geq m\geq 1$,
\begin{equation}\label{stineq}
P_{m-1}(q)P_{n+1}(q)-P_{m}(q)P_{n}(q)\geq_q 0.
\end{equation}
\end{theo}

\proof  Throughout the proof, we simply write $P_n$ for $P_n(q)$.
Let
$$P'_nP_{m-1}-P_nP'_{m-1}=\sum
\limits_{i=0}^{m+n-1}A_iq^i,$$ where $P_n'$ is the derivative of
$P_n$ with respect to $q$.
 We claim that $A_i\geq 0$ for
any $i$. Invoking the recurrence relation \eqref{rec}, the
coefficient of $q^i$ in $P'_nP_{m-1}$ equals
$$
\begin{array}{rcl}
\lefteqn{\sum\limits_{k=0}^i (i-k+1) T(n,i-k+1)T(m-1,k)}&\\[3pt]
&= \sum\limits_{k=0}^i (i-k+1)
(a_1n+a_2(i-k+1)+a_3)T(n-1,i-k+1) T(m-1,k)\\[3pt]
& \quad +\, \sum\limits_{k=0}^i (i-k+1)
(b_1n+b_2(i-k+1)+b_3)T(n-1,i-k)T(m-1,k).
\end{array}
$$
 Again, based on \eqref{rec},
the coefficient of $q^i$ in $P_nP'_{m-1}$ equals
$$
\begin{array}{rcl}
\lefteqn{\sum\limits_{k=1}^{i+1} k T(n,i-k+1)T(m-1,k)}\\[3pt]
&= \sum\limits_{k=1}^{i+1} k
(a_1n+a_2(i-k+1)+a_3)T(n-1,i-k+1)T(m-1,k)\\[3pt]
& \quad +\, \sum
\limits_{k=0}^{i}k(b_1n+b_2(i-k+1)+b_3)T(n-1,i-k)T(m-1,k).
\end{array}
$$

Let
\begin{eqnarray*}
c_k & = & (i-2k+1)(a_1n+a_2(i-k+1)+a_3), \\
d_k & = & (i-2k+1)(b_1n+b_2(i-k+1)+b_3).
\end{eqnarray*}

For $0\leq k\leq \left\lfloor i/2 \right\rfloor$, we find that
$$
\begin{array}{rcl}
c_k+c_{i-k+1} & = & (i-2k+1)(a_1n+a_2(i-k+1)+a_3)\\[5pt]
& & \quad +\,(2k-i-1)(a_1n+a_2k+a_3)\\[5pt]
&=& a_2(i-2k+1)^2,
\end{array}
$$
which is nonnegative. Moreover,
$$
\begin{array}{rcl}
d_k+d_{i-k}&=&(i-2k+1)(b_1n+b_2(i-k+1)+b_3)\\[5pt]
& & \quad +\,(2k-i+1)(b_1n+b_2(k+1)+b_3)\\[5pt]
& = & b_2(i-2k+1)(i-2k)+2(b_1n+b_2(k+1)+b_3),
\end{array}
$$
which is also nonnegative.

 By Lemma \ref{change}, for $0\leq
k\leq \left\lfloor i/2\right\rfloor$, we obtain
\begin{eqnarray}
T(m-1,0)T(n-1,i+1) & \geq & T(m-1,i+1)T(n-1,0),\label{eq-init}\\[5pt]
T(m-1,k)T(n-1,i-k+1) & \geq & T(m-1,i-k+1)T(n-1,k),\quad \label{ineq1}\\[5pt]
T(m-1,k)T(n-1,i-k) & \geq & T(m-1,i-k)T(n-1,k).\label{ineq2}
\end{eqnarray}

We now need to consider the parity of $i$. First, consider the case
when $i$ is odd. Suppose that $i=2j+1$ for some $j$. Clearly, we
have $c_{j+1}=0$. Since all $T(n,k)$ are nonnegative, we get
\allowdisplaybreaks
\begin{align*}
A_i  = & (i+1)(a_1n+a_2(i+1)+a_3)T(m-1,0)T(n-1,i+1)\\
    &    \quad -\,(i+1)(a_1n+a_3)T(m-1,i+1)T(n-1,0)\\
    &    \quad +\, \sum\limits_{k=1}^{i}c_kT(m-1,k)T(n-1,i-k+1)\\
    &    \quad +\, \sum\limits_{k=0}^{i}d_kT(m-1,k)T(n-1,i-k)\\
    &\geq  \sum\limits_{k=1}^{i}c_kT(m-1,k)T(n-1,i-k+1) & \hfill{(\mbox{by \eqref{eq-init}})}\\
    &    \quad +\, \sum\limits_{k=0}^{i}d_kT(m-1,k)T(n-1,i-k)\\
    & \geq  \sum
\limits_{k=1}^{j}(c_k+c_{i-k+1})T(m-1,k)T(n-1,i-k+1)\quad\quad &\hfill{\rule{28pt}{0pt}(\mbox{by \eqref{ineq1}})}\\
&  \quad +\,\sum \limits_{k=0}^{j}(d_k+d_{i-k})T(m-1,k)T(n-1,i-k),
&\hfill{(\mbox{by \eqref{ineq2}})}
\end{align*}
which is nonnegative, since   $c_k+c_{i-k+1}\geq 0$ and
$d_k+d_{i-k}\geq 0$.

The case when $i$ is even can be dealt with via a similar argument.
Suppose that $i=2j$ for some $j$. In this case, we have $d_{j}\geq
0$. Therefore \allowdisplaybreaks
\begin{align*}
A_i    \geq & \sum\limits_{k=1}^{i}c_kT(m-1,k)T(n-1,i-k+1) & \hfill{\rule{50pt}{0pt}(\mbox{by \eqref{eq-init}})}\\
          & \quad +\, \sum\limits_{k=0}^{i}d_kT(m-1,k)T(n-1,i-k)\\
       & \geq  \sum
\limits_{k=1}^{j}(c_k+c_{i-k+1})T(m-1,k)T(n-1,i-k+1)  & \hfill{(\mbox{by \eqref{ineq1}})}\\
&\quad +\,\sum \limits_{k=0}^{j-1}(d_k+d_{i-k})T(m-1,k)T(n-1,i-k) &  \hfill{(\mbox{by \eqref{ineq2}})}\\
&\quad +\, d_j T(m-1,j)T(n-1,j),
\end{align*}
which is nonnegative, since $c_k+c_{i-k+1}\geq 0$, $d_k+d_{i-k}\geq
0$ and  $d_j\geq 0$.

Combining both cases, we are led to the assertion that  $A_i\geq 0$,
namely, for any $n\geq m\geq 1$,
\begin{equation}\label{eq-i}
P'_nP_{m-1}-P_nP'_{m-1}\geq_q 0.
\end{equation}

In view of the recurrence relation \eqref{rec}, we obtain
\allowdisplaybreaks
\begin{align*}
P_m     = & \sum_{k=0}^m T(m,k) q^k\\
        = & \sum_{k=0}^m (a_1m+a_2k+a_3) T(m-1,k) q^k\\
          & \quad +\,\sum_{k=0}^m (b_1m+b_2k+b_3) T(m-1,k-1) q^k\\
        = & \sum_{k=0}^{m-1} (a_1m+a_2k+a_3) T(m-1,k) q^k  \hfill{\rule{76pt}{0pt}(\mbox{by $T(m-1,m)=0$})}\\
          & \quad +\,\sum_{k=1}^m (b_1m+b_2k+b_3) T(m-1,k-1) q^k  \hfill{\rule{30pt}{0pt}(\mbox{by $T(m-1,-1)=0$})} \\
        = & \sum_{k=0}^{m-1} (a_1m+a_2k+a_3) T(m-1,k) q^k\\
          & \quad +\,\sum_{k=0}^{m-1} (b_1m+b_2(k+1)+b_3)q T(m-1,k) q^k\\
        = & \sum_{k=0}^{m-1} (a_1m+a_3) T(m-1,k) q^k+\sum_{k=0}^{m-1} (b_1m+b_2+b_3)q T(m-1,k) q^{k}\\
          & \quad +\,\sum_{k=0}^{m-1} (a_2+b_2q)k T(m-1,k) q^{k}\\
        = & (a_1m+a_3+b_1mq+b_2q+b_3q)P_{m-1}+\sum_{k=0}^{m-1} (a_2q+b_2q^2)k T(m-1,k) q^{k-1}\\
        = &
       (a_1m+a_3+b_1mq+b_2q+b_3q)P_{m-1}+(a_2+b_2q)qP'_{m-1},
\end{align*}
and hence
\[
P_{n+1}=(a_1(n+1)+a_3+b_1(n+1)q+b_2q+b_3q)P_{n}+(a_2+b_2q)qP'_{n}.
\]
Substituting $P_m$ and $P_{n+1}$ into \eqref{stineq} yields
$$
\begin{array}{rcl}
P_{m-1}P_{n+1}-P_{m}P_{n} &=& (a_1+b_1q)(n-m+1)P_{m-1}P_{n}\\[10pt]
& & \quad + \, q(a_2+b_2q)(P'_nP_{m-1}-P_nP'_{m-1}).
\end{array}
$$
Since $a_1,a_2,b_1,b_2\geq 0$, the strong $q$-log-convexity of
$\{P_n(q)\}_{n\geq 0}$ immediately follows from \eqref{eq-i}.\qed

\section{Applications}

In this section, we use Theorem \ref{strong} to show that the Bell
polynomials, the Bessel polynomials, the Ramanujan polynomials and
the Dowling polynomials are strongly $q$-log-convex.

\subsection{The Bell polynomials}

The Bell polynomials  \cite{bell} are defined by
$$B_n(q)=\sum_{k=0}^n S(n,k)q^k,$$
where $S(n,k)$ is the Stirling number of the second kind satisfying
the following recurrence relation:
\[
S(n,k)=kS(n-1,k)+S(n-1,k-1), \quad n, k\geq 1
\]
with $S(0,0)=1$.

\begin{coro}\label{bell}
The Bell polynomials are strongly $q$-log-convex.
\end{coro}

Moreover, one can check that the following polynomials introduced by
Tanny \cite{tanny} are also strongly $q$-log-convex:
$$F_n(q)=\sum_{k=0}^n k!S(n,k)q^k.$$

\subsection{The Bessel polynomials}

The Bessel polynomials are defined by
\[
y_n(x)=\sum
\limits_{k=0}^{n}\frac{(n+k)!}{(n-k)!k!}\left(\frac{x}{2}\right)^k,
\]
and they have been extensively studied; See, for example, Burchnall
\cite{bur1957}, Carlitz \cite{car1957}, and Grosswald
\cite{gro1951}. The Bessel polynomials satisfy the following
recurrence relation \cite{kr1949}:
\begin{eqnarray}
y_n & = & (nx+1)y_{n-1}+x^2y'_{n-1}. \label{recu3}
\end{eqnarray}
Let
\[
T(n,k)=\frac{(n+k)!}{(n-k)!k!}.
\]
From the recurrence (\ref{recu3}), we deduce that
\[
T(n,k)=T(n-1,k)+(2n+2k-2)T(n-1,k-1), \quad n\geq k\geq 1.
\]

\begin{coro}\label{bessel}
The Bessel polynomials are strongly $q$-log-convex.
\end{coro}

\subsection{The Ramanujan polynomials}

The Ramanujan polynomials $R_n(x)$ are defined by the following
recurrence relation:
\begin{equation} \label{id-ra}
R_1(x)=1,\quad R_{n+1}(x)=n(1+x)R_{n}(x)+x^2R_{n}'(x),
\end{equation}
where $R_{n}'(x)$ is the derivative of $R_n(x)$ with respect to $x$.
These polynomials  are related to a  refinement of  Cayley's theorem
due to Shor \cite{Shor}. The connection between the Ramanujan
polynomials and Shor's refinement of Cayley's formula was observed
by Zeng \cite{zeng1999}.  Let $r(n,k)$ be the number of rooted
labeled trees on $n$ vertices with $k$ improper edges. Shor
\cite{Shor}  proved that $r(n,k)$ satisfies the following recurrence
relation:
\begin{equation}\label{ram2}
r(n,k)=(n-1)r(n-1,k)+(n+k-2)r(n-1,k-1).
\end{equation}
where $r(1,0)=1, n\geq 1, k\leq n-1$, and $r(n,k)=0$ otherwise. It
can be seen from \eqref{id-ra} and \eqref{ram2} that $R_n(x)$ is
indeed the generating function of $r(n,k)$, namely,
$$R_n(x)=\sum_{k=0}^n r(n,k)x^k.$$
Dumont and Ramamonjisoa \cite{Dum} independently  found the same
combinatorial interpretation for the coefficients of the Ramanujan
polynomial $R_n(x)$.

Let $r'(n,k)=r(n+1,k)$. Then the triangle $\{r'(n,k)\}_{n\geq k\geq
0}$ satisfies the following recurrence relation
\begin{equation}\label{ram3}
r'(n,k)=nr'(n-1,k)+(n+k-1)r'(n-1,k-1).
\end{equation}

\begin{coro}\label{ram}
The Ramanujan polynomials $R_n(q)$ are strongly $q$-log-convex.
\end{coro}

\subsection{The Dowling polynomials}

The Dowling polynomials are defined as the generating functions of
Whitney numbers of the second kind of Dowling lattices; See
Benoumhani \cite{Ben2}.  As a generalization of the partition
lattice, Dowling \cite{dowling} introduced  a class of geometric
lattices based on finite groups, called the Dowling lattice. Given a
finite group $G$ of order $m\geq 1$, let $Q_n(G)$ be the Dowling
lattice of rank $n$ associated to $G$, and, for $0\leq k\leq n$, let
$W_m(n, k)$ be the Whitney numbers of the second kind of $Q_n(G)$.
The Dowling polynomial $D_m(n;q)$ is defined by
\begin{equation*}
D_m(n;q)=\sum_{k=0}^{n} W_m(n, k)q^k.
\end{equation*}
Benoumhani \cite{Ben2} also introduced the following generalized
polynomials
\begin{equation*}
F_{m,1}(n; q)=\sum_{k=0}^{n} k!\, W_m(n, k)m^k q^k,\quad F_{m,2}(n;
q)=\sum_{k=0}^{n} k!\, W_m(n, k)q^k.
\end{equation*}

Dowling \cite{dowling} proved that the Whitney numbers $W_m(n, k)$
satisfy the following recurrence relation
\begin{equation}\label{eq-wt}
W_m(n, k)=(1+mk)W_m(n-1, k)+W_m(n-1, k-1), \quad n\geq k\geq 1,
\end{equation}
with the boundary conditions:
$$
\begin{array}{rcl}
W_m(n, n) & = & W_m(n, 0)=1, \mbox{ for } n\geq 0,\\[5pt]
W_m(n, k) & = & 0, \mbox{ if } k>n.
\end{array}
$$
Note that, for $m=1$, the Whitney numbers $W_m(n, k)$ of the second
kind are the Stirling numbers $S(n+1,k+1)$ of the second kind; See
\cite{dowling, Ben2}.

From \eqref{eq-wt}, Benoumhani \cite{Ben2} derived that
\begin{equation}\label{eq-fw}
F_{m,2}(n;q) = (q+1)F_{m,2}(n-1;q)+ q(q +m)F'_{m,2}(n-1;q),
\end{equation}
where $F'$ is the derivative of $F$ with respect to $q$. Let
$T(n,k)=k!\, W_m(n, k)$. Then  \eqref{eq-fw} implies that $T(n,k)$
satisfies the following  recurrence relation
\begin{equation}\label{eq-fw1}
T(n, k)=(1+mk)T(n-1, k)+kT(n-1, k-1), \quad n\geq k\geq 1.
\end{equation}

\begin{coro}\label{dm}
The Dowling polynomials $D_m(n;q), F_{m,1}(n; q)$ and $F_{m,2}(n;
q)$ are strongly $q$-log-convex.
\end{coro}

\section{The Bessel transformation}

The objective of this section is to show that the Bessel
transformation preserves log-convexity. The Bessel transformation is
a linear transformation associated with the Bessel polynomials,
which transforms a  sequence $\{z_n\}_{n\geq 0}$ of nonnegative real
numbers into a sequence $\{w_n\}_{n\geq 0}$ given by
$$w_n=\sum_{k=0}^n \frac{(n+k)!}{(n-k)!k!} z_k.$$

For a triangular array $\{T(n,k)\}_{n\geq k\geq 0}$, let $P_n(q)$ be
the row generating function of $T(n,k)$. Given any $n\geq 1$, Liu
and Wang defined the function $\alpha_k(n, i)$ for $0\leq i\leq 2n$
and $0\leq k \leq \lfloor{i\over 2}\rfloor$. If  $0\leq
k<\frac{i}{2}$, we have
\begin{eqnarray*}
\lefteqn{\alpha_k(n,i)=T(n-1,k)T(n+1,i-k)}&&\\[5pt]
&&\hspace{16mm}+T(n+1,k)T(n-1,i-k)-2T(n,k)T(n,i-k).
\end{eqnarray*}
If $i$ is even and $k=\frac{i}{2}$, then we have
\[
\alpha_k(n,i)=T(n-1,k)T(n+1,k)-T(n,k)^2.
\]
Liu and Wang \cite{liuwan2006} found the following connection
between the $q$-log-convexity of $P_n(q)$ and the log-convexity
preserving property of $T(n,k)$.

\begin{theo}[\cite{liuwan2006}]\label{Liuwang2}
Suppose that the triangle $\{T(n,k)\}_{n\geq k\geq 0}$ of positive
real numbers satisfies the following two conditions:

(C1) The sequence of polynomials $\{P_n(q)\}_{n\geq 0}$ is
$q$-log-convex.

(C2) There exists an integer $k'$ depending on $n$ and $i$ such that
$\alpha_k(n,i)\geq 0$ for $k\leq k'$ and $a_k(n,i)<0$ for $k>k'$.

Then the linear transformation $w_n=\sum \limits_{k=0}^{n}
T(n,k)z_k$ preserves  log-convexity.
\end{theo}

We will use the above theorem to prove that the Bessel
transformation preserves log-convexity. For any $n\geq 1$ and $0\leq
i\leq 2n$, we introduce the following polynomials in $x$:
$$
\begin{array}{rcl}
f_1(x) & = & (n+x+1)(n-i+x)(n+x)(n-i+x+1),\\[5pt]
f_2(x) & = & (n-x)(n+i-x+1)(n-x+1)(n+i-x),\\[5pt]
f_3(x)  & = & (n-x+1)(n+i-x)(n+x)(n-i+x+1).
\end{array}
$$
Let \begin{equation}\label{fx}
 f(x)=f_1(x)+f_2(x)-2f_3(x).
\end{equation}

\begin{lemm}\label{lemm-con}
For any $n\geq 1$, $0\leq i\leq 2n$ and $0\leq k\leq \lfloor
\frac{i}{2}\rfloor$, let
\begin{eqnarray*}
\beta_k(n,i) & = & T(n+1,k)T(n-1,i-k)+T(n+1,i-k)T(n-1,k)\\
&&\quad -\,2T(n,i-k)T(n,k),
\end{eqnarray*}
where
\begin{equation}\label{eq-b}
T(n,k)=\frac{(n+k)!}{(n-k)!k!}.
\end{equation}
Then, for given $n$ and $i$, there exists an integer $k'$ depending
on $n$ and $i$ such that $\beta_k(n,i)\geq 0$ for $k\leq k'$ and
$\beta_k(n,i)\leq 0$ for $k> k'$.
\end{lemm}

\proof Suppose that  $n$ and $i$ are given. Clearly, if $k< i-n-1$,
then $n< (i-k)-1$ and $\beta_k(n,i)=0$. If $k=i-n-1$, then
\[ \beta_{k}(n,i)=T(n+1,n+1)T(n-1,i-n-1)\geq 0.\]
 Therefore, it suffices to determine the sign
of $\beta_k(n,i)$ for $i-n-1< k\leq \lfloor \frac{i}{2}\rfloor$.

By \eqref{eq-b},  we have
$$
\begin{array}{rcl} \beta_k(n,i) & = &
\displaystyle{\frac{(n+1+k)!}{(n+1-k)!\,k!}}\times
\frac{(n-1+i-k)!}{(n-1-i+k)!\,(i-k)!}\\[12pt]
& & \quad  +\,
\displaystyle{\frac{(n+1+i-k)!}{(n+1-i+k)!\,(i-k)!}\times \frac{(n-1+k)!}{(n-1-k)!\,k!}}\\[12pt]
&& \quad -\,\displaystyle{\frac{2(n+i-k)!}{(n-i+k)!\,(i-k)!}\times \frac{(n+k)!}{(n-k)!\,k!}}\\[12pt]
&=&
\left(\displaystyle{\frac{(n+k+1)(n-i+k)}{(n-k+1)(n+i-k)}+\frac{(n-k)(n+i-k+1)}{(n+k)(n-i+k+1)}-2}\right)\\[12pt]
&& \quad \times\, \displaystyle{\frac{(n+i-k)!}{(n-i+k)!\,(i-k)!}}\times \displaystyle{\frac{(n+k)!}{(n-k)!\,k!}}\\[12pt]
&=&\displaystyle{\frac{f(k)T(n,i-k)T(n,k)}{f_3(k)}}.
\end{array}
$$

Note that for $n\geq 1$, $0\leq i\leq 2n$ and $i-n-1< k\leq \lfloor
\frac{i}{2}\rfloor$, we have
$$f_3(k)>0, \quad T(n,i-k)T(n,k)\geq 0.$$

By the definition (\ref{fx}) of $f(x)$, we find
$$f'(x)=2(2x-i)g(x),$$
where
$$g(x)=2(2+8n^2-i+8n).$$
Thus, for $i\leq 2n$ and $x\leq \frac{i}{2}$, we have
$$g(x)\geq 2(2+8n^2+6n)>0, \quad f'(x)\leq 0.$$
Therefore, $f(x)$ is decreasing in $x$ on the interval
$(-\infty,\frac{i}{2}]$. This implies that there exists an integer
$k'$ such that $\beta_k(n,i)\geq 0$ for $k\leq k'$ and
$\beta_k(n,i)\leq 0$ for $k> k'$. \qed

Combining Corollary \ref{bessel}, Lemma \ref{lemm-con} and Theorem
\ref{Liuwang2}, we deduce the following theorem.

\begin{theo}
If $\{z_k\}_{k\geq 0}$ is a log-convex sequence of positive real
numbers, then the sequence $\{w_n\}_{n\geq 0}$ defined by
$$w_n=\sum_{k=0}^n
\frac{(n+k)!}{(n-k)!k!} z_k$$ is also log-convex.
\end{theo}

\noindent {\bf Acknowledgments.} This work was supported by  the 973
Project, the PCSIRT Project of the Ministry of Education, the
Ministry of Science and Technology, and the National Science
Foundation of China.

\end{document}